# Volterra model-based control for nonlinear systems via Carleman linearization

Dhruvi Bhatt*

*Department of Electrical Engineering, Sardar Vallabhbhai National Institute of Technology, Surat, Gujarat 395007, India. E-mail address: [shn411@gmail.com](shn411@gmail.com)*

Shambhu Nath Sharma

*Department of Electrical Engineering, Sardar Vallabhbhai National Institute of Technology, Surat, Gujarat 395007, India. E-mail address: [snsvolterra@gmail.com](snsvolterra@gmail.com)*

**Abstract**

   This paper presents detailed insights of embedding Carleman linearization into nonlinear systems for designing Volterra model-based control technique. Volterra series method is a competent mathematical tool, which extends the convolution integral for linear systems to nonlinear systems. First, we utilize the Carleman linearization technique to arrive at the bilinear approximation of the nonlinear system. Secondly, the third-order Volterra model representation is computed from the Carleman bilinearized model. Then, Volterra model-based control strategy is developed. The proposed method is effectuated on the benchmark van de Vusse reactor – exhibiting non-minimum phase response. The open and closed-loop simulation results are presented, demonstrating the superiority and practical utility of the proposed method.

## 1. Introduction

   In real-world, physical systems are inherently nonlinear in behaviour and nature. Nonlinear control problems have drawn immense interest to the people of engineering discipline. Foremost reasons behind increasing interest in nonlinear control systems incorporate improvement in linear control systems, in analysis of nonlinearities and simplicity in designing. More generally, a nonlinear system can be represented by a Volterra model, which is based on the Volterra series (Burt & De Morais Goulart, 2018). Herein, for the analysis of nonlinear dynamical system the Volterra series approach is adopted. Volterra and Wiener series are two classes of polynomial representations of nonlinear systems. They are perhaps the best and most widely used nonlinear system representation technique for highly non-linear systems. This approach to the systematic characterization of non-linear system is given by a great mathematician Vito Volterra back in

1887. The Volterra series representation is an extension of the classical linear system representation and it represents an important model for the representation, analysis, and synthesis of non-linear systems.

Nonlinear problems have drawn great interest and extensive attention from engineers, physicists and mathematicians and many other scientists because most real systems are inherently nonlinear in nature. To model and analyze nonlinear systems, many mathematical theories and methods have been developed, including Volterra series (Cheng et al., 2017; George, 1959; Ikehara, 1951). In many of the problems encountered in nonlinear systems analysis, it is difficult to find explicit closed-form expressions for the solution. Volterra series is one of the most widely used and well-established methods. The intention behind using the Volterra functional series is fundamentally to evaluate solutions which cannot be obtained in closed form, and therefore the solution obtained from this approach will involve an approximation.

The theory of differential equations is a fundamental tool in physics, engineering and other mathematically based sciences. Many natural laws and models of natural phenomena are described by nonlinear systems of differential equations. Many nonlinear problems are difficult or impossible to solve in closed form and therefore the construction of such solutions is nontrivial. Carleman developed a technique to linearize nonlinear systems of equations that could lead to analytical solutions of nonlinear problems. Several techniques in the area of nonlinear model reduction have been investigated (Burt & De Morais Goulart, 2018). One such technique is the one based on the construction of tractable bilinear models of nonlinear systems. This is motivated by the fact that bilinear systems are ubiquitous among nonlinear systems with linear inputs.

Here, we describe a partial historic account of contributions towards the control of nonlinear systems. It can be traced back to the work of the Italian mathematician Vito Volterra (1959) about theory of analytic functional in 1887. Norbert Wiener (1942) applied his theory of Brownian motion to investigate the integration of Volterra analytic functional and firstly used it for system analysis. Volterra suggested a series representation which provides an explicit description of the input-output relationship for the class of non-linear sym. Torsten Carleman (1932) showed a finite dimensional system of nonlinear differential equation can be embedded into an infinite system of linear differential equations. This method is named as Carleman Linearization. Martin Schetzen (1980) presented a theory that wider class of non-linear dynamical system can be represented by a Volterra model, which is based on Volterra series representation. There are wide applications of Volterra series in different disciplines from aero-elastic system to mechanical system. Silva (2005) showed an application of Volterra theory for identification of aero-elastic system. Similarly in electrical engineering discipline Ronnow et al. (2007) demonstrated characterization of nonlinear memory effects in radio-frequency power amplifiers. Korenberg and Hunter (1996) showed an application of Volterra theory and identification of kernals on biological systems. For mechanical systems, Jing and Lang (2009) applied Volterra series theory on damping system subjected to harmonic input. Khan and Vyas (2001) explained an application of Volterra theories on rotor bearing system. Joaquín Collado (2008) showed the Carleman linearization technique representation for which the design of controllers can be performed directly.

This paper presents detailed analysis of IMC based Volterra model control of non-linear system. Volterra Series method is a competent mathematical tool which extends the convolution integral for linear systems to non-linear systems. The original nonlinear system is replaced by a bilinear system via formulation of Carleman linearization technique. This technique is used for identification of Volterra kernals which surpass the computational hurdles. Then, the Volterra model is computed from the approximated Carleman bilinearized system. In this paper, Carleman linearization is applied to a benchmark control problem of van de Vusse reactor which incorporates square nonlinearity. In this case study, third-order Volterra model is obtained and it is compared with second-order and first-order Volterra models. Realization of first, second and third-order Volterra model simulation studies for both open-loop as well as closed-loop is carried out. Control investigation is presented by designing an IMC-Volterra controller scheme for van de Vusse reactor. The illustrative simulation results indicate the superiority of third-order Volterra model-based controller which shows closeness to the true non-linear model and ideal system.

## 2. Background

### 2.1 Volterra series model

The Volterra series method is a type of functional series which relates the system input, to the system output (Schetzen 1980). The Volterra series gives a mathematical representation of the solution in the form of an expanding infinite series of integrals which encompasses the nonlinearities of the system. The response of a causal, time-invariant, finite-memory, nonlinear system to an arbitrary input may be expressed as an infinite sum of multidimensional convolution integrals of increasing order.

$$y(t) = y_0 + \int_0^t h_1(\tau_1)u(t-\tau_1)d\tau + \int_0^t \int_0^t h_2(\tau_1,\tau_2) u(t-\tau_1)u(t-\tau_2)d\tau_1,\tau_2$$

$$+ \sum_{i=3}^{\infty} \int_0^t \cdots \int_0^t h_i(\tau_1,\cdots,\tau_i) u(t-\tau) \cdots u(t-\tau_i)d\tau_1 \cdots d\tau_i,$$

where $y(t)$ is an $m \times 1$ response vector, $u(t)$ is an $r \times 1$ input vector, $h_1(\tau_1)$ is a first-order kernel, $h_2(\tau_1,\tau_2)$ is a second-order kernel, etc.

Representations of the responses of nonlinear systems having the form given in equation are called Volterra series and the kernels appearing in them are called Volterra kernels. The first integral appearing in equation is the linear convolution integral that gives the linear contribution to the total response. The second integral in the series may be regarded as a second order convolution integral that provides the second-order contribution to the total response. Similar interpretations may be made for the successively higher-order integrals in the series. The first-order kernel $h_1(\tau_1)$ characterizes the linear behavior of the system and is identical to the unit impulse response of the linear (or linearized) system. The second-order kernel $h_2(\tau_1, \tau_2)$ characterizes the behavior of the nonlinear system to two separate unit impulses applied at two varying points in time.

## 2.2 Carleman linearization

Carleman linearization is a widely applied method for approximating a linear-analytic system by a bilinear system having exactly the same associated Volterra model up to a prescribed order. In fact, because of this property, Carleman linearization is also used as the first step in obtaining the Volterra model of a linear-analytic system (Ekman 2005, Burt 2018).

Carleman linearization technique has been developed to transform sets of polynomial ordinary differential equations into infinite dimensional linear system representation. After choosing suitable finite orders of Carleman linearization finite dimensional linear and bilinear system representation are obtained for which the design of controllers can be performed straight forward. Carleman linearization is applied to extend the state vectors of original state variables. An advantage of this method is the embedding of dynamics of nonlinear systems in corresponding bilinear representations. The Carleman Linearization was proposed for autonomous nonlinear analytic systems almost 80 years ago (Carleman 1932).

One approach for calculating the Volterra kernels in an input/output representation is the so called Carleman linearization. Sometimes the approach is simply called the bilinearization method. The Carleman linearization method is used for computing kernels for nonlinear state equations which are affine in the input.

$$\dot{x} = f(x,t) + g(x,t)u, \ x(0) = x_0. \tag{1}$$

The input signal u is scalar and it is assumed that (1) is a linear analytic state equation i.e. $f(x,t)$ and $g(x,t)$ are analytic in x and continuous in t. Then, it can be shown that there exists a Volterra system representation for (1), that converges on $t \in [0,T]$ for sufficiently small inputs signal (Rugh 1980). For simplicity we will consider a scalar output $y(t) = Cx(t)$.

*Assumption:* Without loss of generality, we assume that $x_0$ is at the origin and $y(x_0)$ is equal to zero as well (Hashemian 2015).

For a given $A \in R^{n \times m}$ and $B \in R^{s \times t}$, their Kronecker product (Brewer 1978) can be written as

$$A \otimes B = \begin{pmatrix} a_{11}B & a_{12}B & \cdots & a_{1m}B \\ a_{21}B & a_{22}B & \cdots & a_{2m}B \\ \vdots & \vdots & \ddots & \vdots \\ a_{n1}B & \cdots & \cdots & a_{nm}B \end{pmatrix}.$$

To achieve the carleman linearized bilinear model we replace the right side of the state equation by adopting the Kronecker-product notation, we can write

$$\frac{d}{dt}x^\otimes = \begin{pmatrix} A_{11} & A_{12} & \cdots & A_{1N} \\ 0 & A_{21} & \cdots & A_{2N-1} \\ \vdots & \vdots & \cdots & \vdots \\ 0 & 0 & \cdots & A_{N1} \end{pmatrix} x^\otimes + \begin{pmatrix} B_{11} & B_{12} & \cdots & A_{1N-1} & 0 \\ 0 & B_{21} & \cdots & B_{2N-2} & 0 \\ \vdots & \vdots & \cdots & \vdots \\ 0 & 0 & \cdots & B_{N0} & 0 \end{pmatrix} x^\otimes u + \begin{pmatrix} B_{10} \\ 0 \\ \vdots \\ 0 \end{pmatrix} u. \quad (2)$$

$$y(t) = (c(t) \ 0 \ \cdots \ 0)x^\otimes + y_0(t), \quad (3)$$

where, $x^\otimes = (x^{(1)} \ x^{(2)} \ \cdots \ x^{(N)})^T$. Suppose the vector $x = (x_1 \ x_2 \ \cdots \ x_n)^T$, then $x \otimes x$ becomes

$$(\underbrace{x_1^2, x_1x_2, x_1x_3, \cdots x_1x_n}, \underbrace{x_2x_1, x_2^2, x_2x_3, \cdots x_2x_n}, \underbrace{x_3x_1, x_3x_2, x_3^2, \cdots x_3x_n}, \cdots \underbrace{x_nx_1, x_nx_2, \cdots x_n^2}).$$

## 2.3 Volterra kernel representation and realization

In this section, we present the Volterra kernels in frequency domain setting. The transformation from time-domain Volterra kernel to frequency-domain kernels involves two-dimensional Laplace transform, which brings unnecessary mathematical intricacies. Thus, without dwelling into greater details, we utilize the results (4)-(6) to achieve the Volterra kernels in frequency-domain (Sain *et al.*, 1990).

$$H_1(s) = c^T(sI_1 - A_{11})^{-1}b_1, \quad (4)$$

$$H_2(s_1, s_2) = c^T((s_1 + s_2)I_2 - A)^{-1}N(s_1 - A)^{-1}b, \quad (5)$$

$$H_3(s_1, s_2, s_3) = c^T((s_1 + s_2 + s_3)I_3 - A)^{-1}N((s_1 + s_2)I_2 - A)^{-1}N(sI_1 - A)^{-1}b, \quad (6)$$

In general, from the induction principle, we get

$$H_n(s_1, s_2, \ldots, s_n) = C^T(((s_1 + s_2 + \cdots + s_n)I_n - A)^{-1}N((s_1 + s_2 + \ldots + s_{n-1})I_{n-1}$$
$$-A)^{-1}N \cdots N((s_1I - A)^{-1}B)) \quad (7)$$

Combing (4) with (2) and (3), we get the first-order portion of the Volterra series representation of the Carleman bilinearized ODEs. Similarly, the second-order and third-order portion can be obtained by (5) and (6), respectively. For nth order portion (7) can be utilized, however, we restrict our discussion to third-order portion.

Realization of Volterra kernels for nonlinear systems is a challenging task. Thanks to Schetzen (1965), for the realization of $n$th order portion of the volterra model $H_n(s_1, s_2, \ldots, s_n)$, we express the above Eq. (7) as a ratio of two polynomials (8). The synthesis of nth-order Volterra kernels is shown in Fig. 1.

$$H_n(s_1, s_2, \ldots, s_n) = \frac{P_1(s_1)P_2(s_2)P_3(s_1+s_2)\cdots P_n(s_1+s_2+\cdots s_n)}{Q_1(s_1)Q_2(s_2)Q_3(s_1+s_2)\cdots Q_n(s_1+s_2+\cdots s_n)}. \tag{8}$$

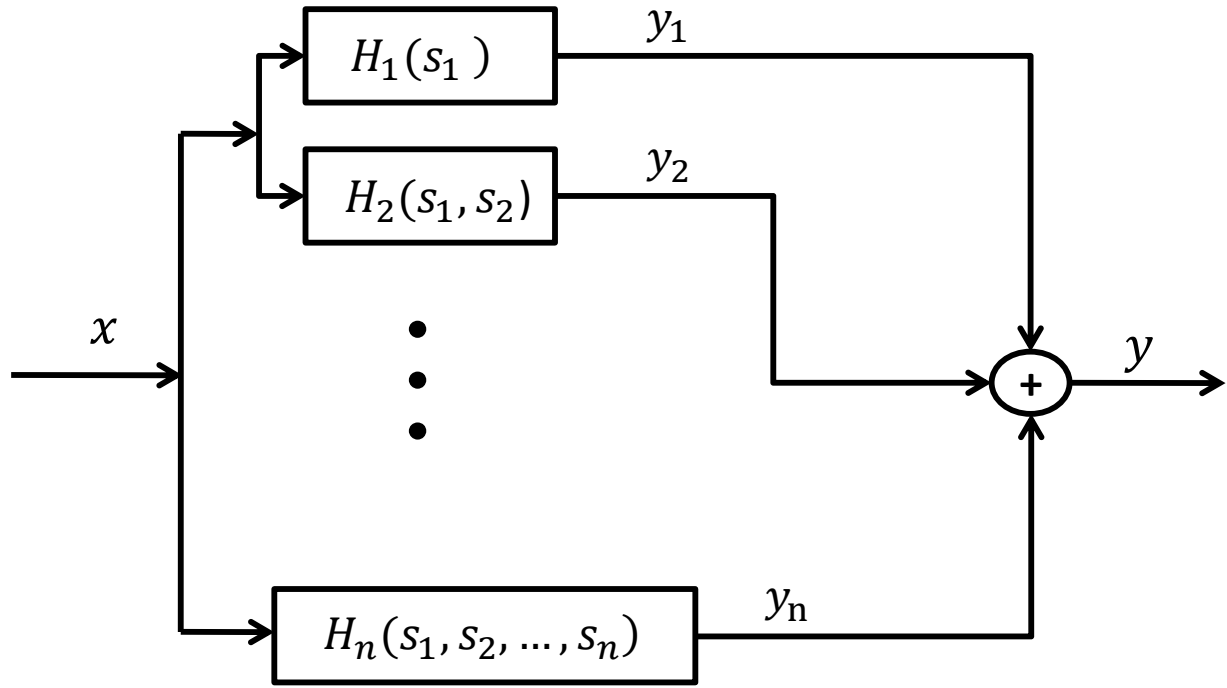

Fig. 1. Synthesis of nth-order Volterra kernels.

For the brevity of presentation, consider case where $n = 2$. Then the second order Volterra model can be expressed as

$$H_2(s_1, s_2) = \frac{P_2(s_1, s_2)}{Q_2(s_1, s_2)}. \tag{9}$$

This second-order Volterra kernel (9) can be realized by a finite number of multipliers. Now $Q_2(s_1, s_2)$ can be represented as a combination of $F_a(s_1), F_b(s_2)$ and $F_c(s_1 + s_2)$, i.e.

$$Q_2(s_1, s_2) = F_a(s_1)F_b(s_2)F_c(s_1 + s_2). \tag{10}$$

The final step is to expand $H_2(s_1, s_2) = \frac{P_2(s_1,s_2)}{F_a(s_1)F_b(s_2)F_c(s_1+s_2)}$ by partial fraction and express as a combination of multipliers.

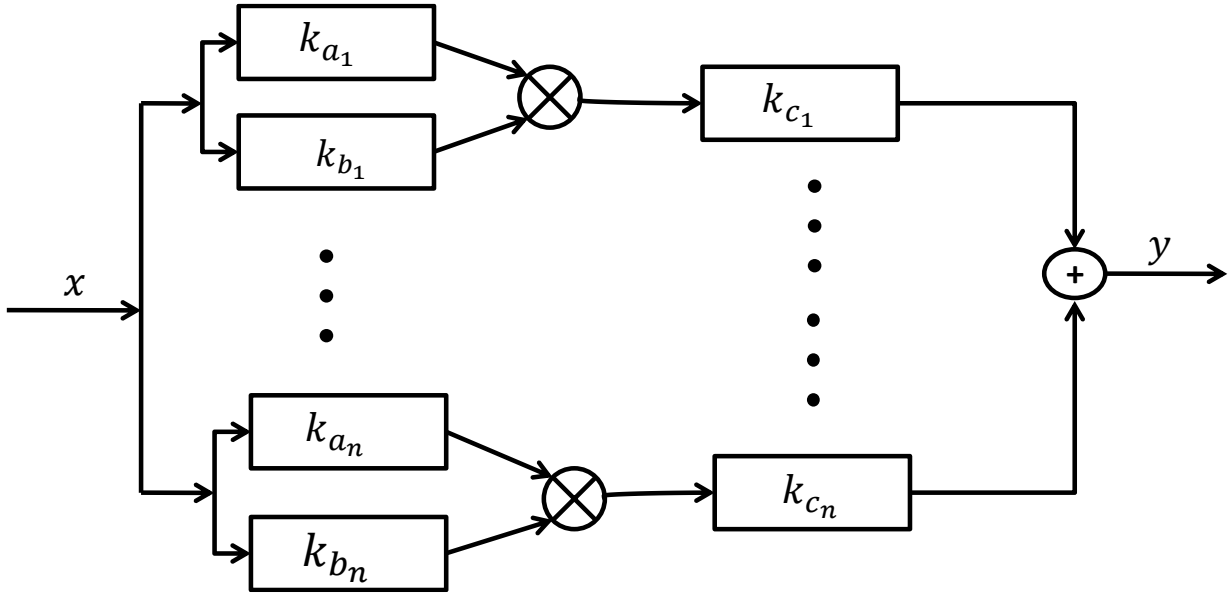

Fig. 1: Realization for second-order Volterra kernels.

Figures (2) and (3) are the schematic block diagrams of realization of second and third-order Volterra kernels, respectively. Here, Figure 2 shows the second order systems with N multipliers. Second-order kernel of the system is denoted by $H\_2(s\_1, s\_2)$, where $k\_(a\_i), k\_(b\_i), k\_(c\_i)$ are the multipliers as a function of $s\_1, s\_2$, and $s\_1 + s\_2$, respectively.

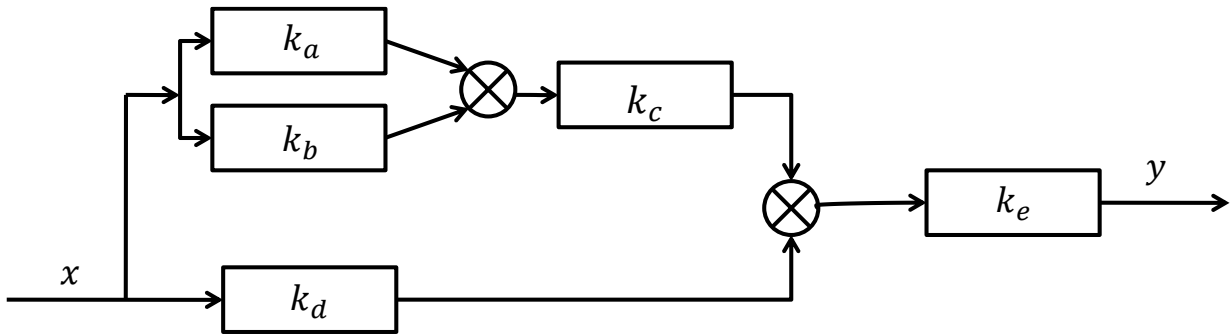

Fig. 2: Realization for Third-order Volterra kernels.

Here, we write steps to synthesize Volterra kernels.

- Express $H_n(s_1, s_2, \ldots, s_n)$ as a ratio of two polynomials.
- Represented the denominator polynomial as a combination of individual polynomials upto $n$th-order.
- Expand $H_n(s_1, s_2, \ldots, s_n)$ by partial fraction and express them as a combination of multipliers.

*Example of Second-order kernel synthesis*

Suppose, we desire to synthesis a kernel $H_2$ having the following structure

$$H_2(s_1, s_2) = \frac{P_2(s_1, s_2)}{Q_2(s_1, s_2)}. \tag{11}$$

$$= \frac{s_1^2 s_2 + 4s_1^2 + 12s_1 s_2 + s_1 s_2^2 + 2s_2^2 + 12s_1 + 48}{s_1^2 s_2 + 4s_1^2 + 12s_1 s_2 + s_1 s_2^2 + 2s_2^2 + 32s_1 + 20s_2 + 48}. \tag{12}$$

This above-mentioned kernel can be synthesis by a combination of finite numbers of multipliers. Thus, following the procedure described above $Q_2(s_1, s_2)$ is expressible in the form of $F_a(s_1), F_b(s_2)$ and $F_c(s_1 + s_2)$. To determine them let take $s_1 = 0$ and $s_2 = s$, we get

$$Q_2(0, s) = F_a(0) F_b(s) F_c(s) = 2s^2 + 20s + 48 = 2(s + 4)(s + 6). \tag{13}$$

Now, let's take $s_1 = s$ and $s_2 = 0$, we get

$$Q_2(s, 0) = F_a(s) F_b(0) F_c(s) = 4(s + 2)(s + 6). \tag{14}$$

Finally, take $s_1 = s$ and $s_2 = -s$,

$$Q_2(s, -s) = F_a(s) F_b(-s) F_c(0) = 6(s + 2)(-s + 4). \tag{15}$$

On comparing the common zeros of the above three eqns. (13)-(15), we get

$$Q_2(s_1, s_2) = (s_1 + 2)(s_2 + 4)(s_1 + s_2 + 6). \tag{16}$$

Embedding the above Eq. (16), into Eq. (12) and expanding by partial fractions, we obtain

$$H_2(s_1, s_2) = \frac{P_2(s_1, s_2)}{Q_2(s_1, s_2)} = 1 - \frac{20(s_1 + s_2)}{(s_1 + 2)(s_2 + 4)(s_1 + s_2 + 6)}. \tag{17}$$

Now the desired second order system can be synthesized as a combination of the two basic second-order systems as shown in Fig. 4.

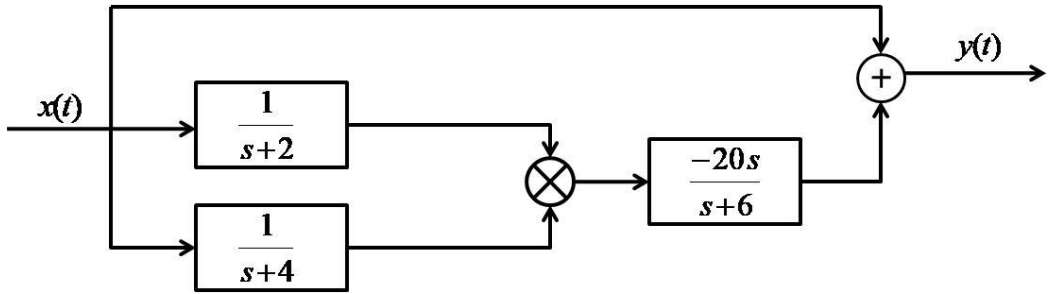

Fig. 3: Realization for second-order Volterra system for specific case.

## 2.4 Internal model control

This section briefly presents Internal Model Control (IMC), and then the IMC-Volterra controller structure for $n$th order Volterra model is sketched. The IMC structure is shown in Fig. 5. The internal model control is a comprehensive model-based controller design method (Garcia & Manfred, 1982)

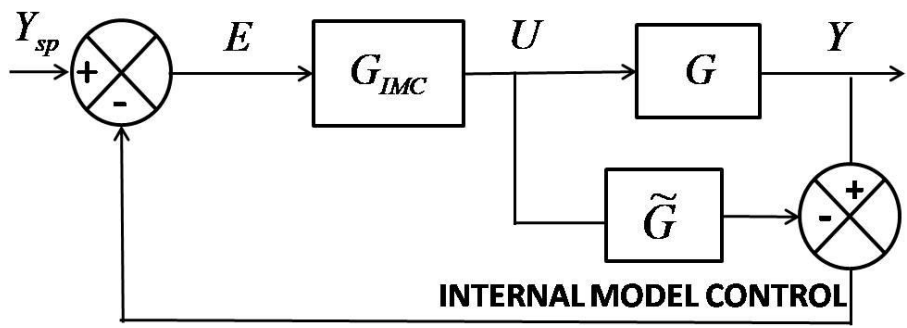

Fig. 5: Linear IMC closed-loop structure

From the block diagram of Fig. (5), the input signal U and the output signal becomes

$$U = G_{IMC}E, Y = GG_{IMC}E \text{ and } E = Y_{sp} - Y^*.$$

where $G_{IMC}(s) = F(s)\tilde{G}_-^{-1}(s)$, $G_c(s) = \dfrac{G_{IMC}(s)}{1 - G_{IMC}(s)\tilde{G}(s)}$. $G_{IMC}$ is the internal model controller consisting of a filter $F(s) = \dfrac{1}{(1+\tau_c s)^r}$ and the inverse of minimum phase portion of the plant model $\tilde{G}$.

## 3. Problem formulation and main result

The IMC model-based design technique leads to a perfect closed-loop performance under the following assumptions. (i) The model is perfect $G = \tilde{G}$. (ii) The controller is equal to the inverse of the model $G_{IMC} = \tilde{G}^{-1}$. (iii) The steady state gain of the controller is equal to that of the model inverse $G_{IMC}(0) = \tilde{G}(0)^{-1}$. However, for IMC controller based on Volterra series model above assumptions does not hold. The case becomes non-trivial and warrants investigations. Hence, we now proceed to develop an explicit expression of the Volterra model-based controller in terms of appropriate generalized inverse. Fig. 6 shows a proposed flowchart for the concerned control technique. First, we embed Carleman linearization in the nonlinear system model to get the bilinearized model. Then, we achieve the Volterra model representation for the Carleman biliniearized model. Finally, the conventional IMC design is reformulated to the IMC-Volterra model-based controller structure.

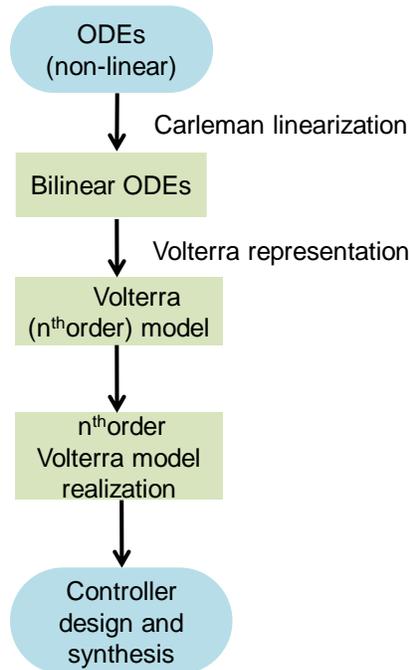

Fig. 6. Flow of the proposed control technique.

We now reformulate the problem as follows. On extending the linear IMC concept to the $n$th order Volterra series model the following block diagrams can be sketched, see Figs. (7)-(9).

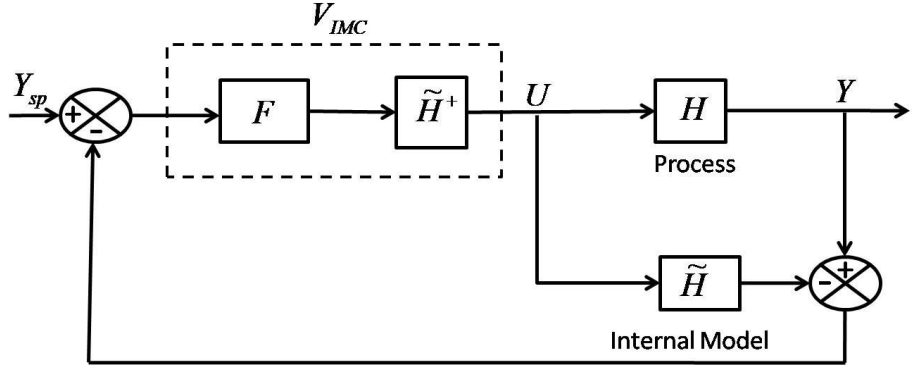

Fig. 7. IMC-Volterra closed-loop structure

From the block diagram, the following relationship follows

$$V_{IMC} = F\widetilde{H}^+ \qquad (18)$$

where $H^+(s_1, s_2, \cdots, s_n)$ is a generalized inverse of the nth order Volterra series. Furthermore, $V_{IMC}(s)$ can be decomposed into a linear model-based controller with a sequence of nonlinear correction terms up to nth order. For the brevity of presentation consider a third order Volterra model, the controller structure $V_{IMC}(s)$ can decomposed into: (i) a linear controller $C(s) = H_1^+(s)$ and (ii) a feedback loop with third-order Volterra model to account for the plant nonlinearity. The alternative implementable scheme for the above is displayed in Fig. (8).

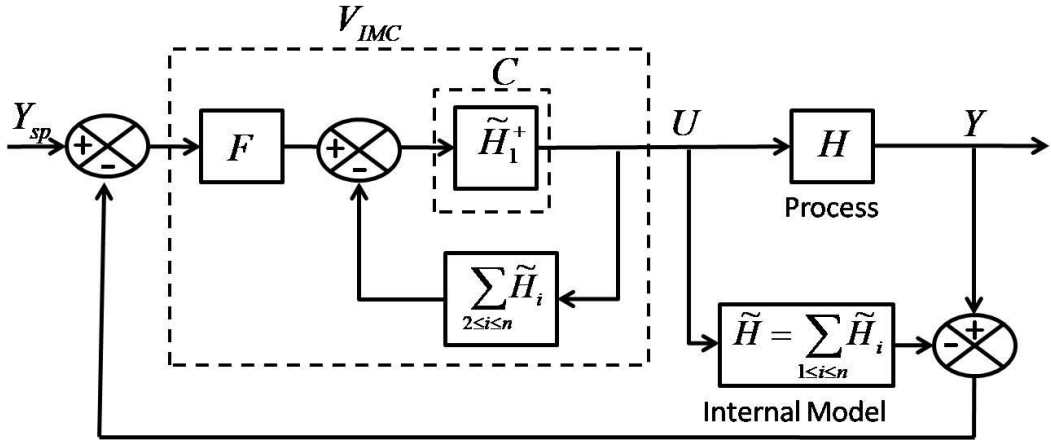

Fig. 8: A closed-loop structure with IMC-Volterra controller

Fig. (8) shows the block diagram of closed-loop system with IMC-Volterra controller for system represented with nth order Volterra model. From the block diagram, $V_{IMC}(s)$ can be recast as

$$V_{IMC} = \frac{F\widetilde{H}_1^+}{1+\widetilde{H}_1^+ \sum_{2\leq i \leq n}\widetilde{H}_i}. \qquad (19)$$

Remark 2: Eq. (18) uses a generalized inverse of the linear model (first order Volterra model). There arise situations in which, the $H_1(s)$ has no inverse or if exist, the true inverse is not possible, for e.g. the linear model consist of inverse response due to the presences of non-minimum phase open loop zeros. Here, the inverse of $\widetilde{H}_1(s)$ is not an exact inverse. Hence, the resulting controller will not be stable or causal. To obtain, a stable linear controller $C, \widetilde{H}_1$ can be factored as a partitioned process with one part as allpass and another as minimum phase contribution. $\widetilde{H}_1$ can be outlined as

$$\widetilde{H}_1 = \widetilde{H}_1^A \widetilde{H}_1^M,$$

where $\widetilde{H}_1^A$ denotes the allpass and $\widetilde{H}_1^M$ is the minimum phase portion. Furthermore, the stable linear controller $\bar{C}$, an pseudo-inverse can be obtained by inverting the minimum phase part of the transfer function $\widetilde{H}_1$, i.e.

$$\bar{C} = F(\widetilde{H}_1^M)^{-1},$$

where $F$ is a first order filter multiplied with the controller $\bar{C}$ to make it realizable. The block diagram in Fig. (9) is a direct consequence of the above change. Fig. (8) and Fig. (9) are equivalent in nature and Fig. (9) is a realizable and implementable version of Fig. (8). The feedback path to the controller $V_{IMC}$ is now modified as $C^* = F^{-1} \sum_{2 \leq i \leq n} \widetilde{H}_i$. Doyle et al., (1995a) showed that the controller which is designed by employing allpass factorization for the partition of the system gives a superior closed-loop response in contrast to the linear controller designed using conventional IMC for non-minimum phase system.

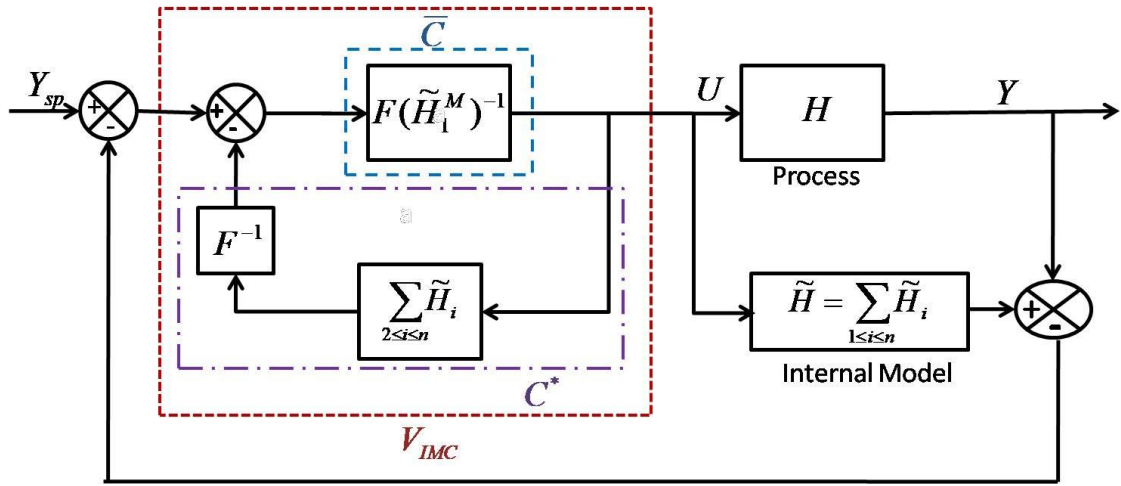

Fig. 9: An alternative realizable closed-loop structure with IMC-Volterra controller

## 4. Case study: van de Vusse reactor

In this section we consider the application of IMC controllers based on Volterra model to a van de Vusse reaction taking place in a continuous stirred tank reactor CSTR. From control system point of view, it is used as a benchmark problem (Dorado *et al.* 2007) for nonlinear processes to demonstrate various process control methodologies.

Consider a CSTR which is operating at a constant temperature (isothermal). The volume is also assumed constant (Doyle *et al.* 1995). The reaction scheme consists of irreversible reactions, see Fig. 1. The feed stream contains only component A. Chemical Reactions of Isothermal CSTR with the kinetics governed by van de Vusse reactions given by,

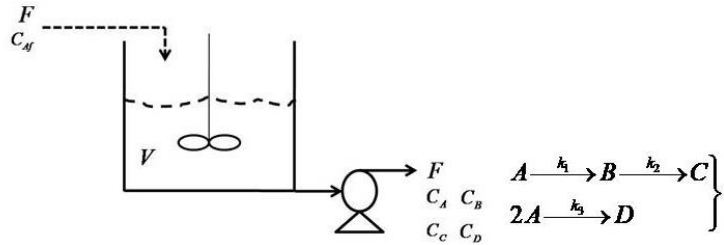

Fig.10: CSTR with van de Vusse reaction scheme

The mathematical model of van de Vusse reactor is represented by the following state equations (Doyle *et al.* 1995)

$$\frac{dC_A}{dt} = \frac{F}{V}(C_{Af} - C_A) - k_1 C_A - k_3 C_A^2. \tag{20}$$

$$\frac{dC_B}{dt} = \frac{FC_B}{V} + k_1 C_A - k_2 C_B, y = C_B. \tag{21}$$

The main objective of presented control problem is to focus on regulating the concentration of component $B$ by manipulating the inlet flow rate. Suppose the concentration of component $A$ is denoted by $x_1$, the concentration of component $B$ is denoted by $x_2$ and the inlet flow rate by $u$. The operating points and model parameter values are listed in Table 1. Substituting the value of parameters in Eq. (4) and Eq. (5) we get,

$$\dot{x}_1 = -50x_1 - 10x_1^2 + u(10 - x_1).$$

$$\dot{x}_2 = 50x_1 - 100x_2 + u(-x_2), y = x_2.$$

Table1: Operating points and model parameters

| Parameters | Value | Parameters | Value |
|---|---|---|---|
| $C_{A0}$ | $3.0 mol l^{-1}$ | $k_2$ | $100\ h^{-1}$ |
| $C_{B0}$ | $1.12 mol l^{-1}$ | $k_3$ | $10 l mol^{-1} h^{-1}$ |
| $\dfrac{F}{V}$ | $34.3 h^{-1}$ | $C_{Af}$ | $10\ mol l^{-1}$ |
| $k_1$ | $50\ h^{-1}$ | $V$ | $1 l$ |

Now on applying the method Carleman linearization to (6) and (7) as described in *Section 1.2*, we get bilinear linear approximation of the nonlinear van de Vusse model, i.e.

$$\frac{dz^{\otimes}}{dt} = \begin{pmatrix} \dfrac{dz_1}{dt} \\ \dfrac{dz_2}{dt} \\ \dfrac{dz_1 z_2}{dt} \\ \dfrac{dz_1^2}{dt} \\ \dfrac{dz_2^2}{dt} \end{pmatrix} = \begin{pmatrix} \dot{z}_1 \\ \dot{z}_2 \\ z_1 \dot{z}_1 + \dot{z}_1 z_2 \\ 2 z_1 \dot{z}_1 \\ 2 z_2 \dot{z}_2 \end{pmatrix},$$

with $z^{\otimes} = (z_1\ z_2\ z_1 z_2\ z_1^2\ z_2^2)^T$, where $\otimes$ denotes Kronecker product. Using the model parameters described in Table 1 and Eq. (8), we arrive at the Carleman linearized bilinear approximation of the van de Vusse reactor model, i.e.

$$\dot{z}_1 = -144.3 z_1 - 10 z_1^2 + 7\tilde{u} - \tilde{u} z_1. \tag{22a}$$

$$\dot{z}_2 = 50 z_1 - 134.3 z_2 - \tilde{u} z_2 - 1.12 \tilde{u}. \tag{22b}$$

$$\dot{z}_1 \dot{z}_2 = 50 z_1^2 - 278.6 z_1 z_2 - 2\tilde{u} z_1 z_2 - 1.12 \tilde{u} z_1 + 7\tilde{u} z_2. \tag{22c}$$

$$\dot{z}_1^2 = -2.88.6 z_1^2 + 14 \tilde{u} z_1 - 2\tilde{u} z_1^2. \tag{22d}$$

$$\dot{z}_2^2 = 100 z_1 z_2 - 268.6 z_2^2 - 2\tilde{u} z_2^2 - 2.24 \tilde{u} z_2. \tag{22e}$$

Rearranging the Eqs. (22a-22e) in bilinear form using Eqs. (2) and (3) we get,

$$\dot{z} = Az + Nz\tilde{u} + b\tilde{u}, y = c^T z. \tag{23}$$

$$\begin{pmatrix} \dot{z}_1 \\ \dot{z}_2 \\ \dot{z}_1\dot{z}_2 \\ \dot{z}_1^2 \\ \dot{z}_2^2 \end{pmatrix} = \begin{pmatrix} -144.3 & 0 & 0 & -10 & 0 \\ 50 & -134.3 & 0 & 0 & 0 \\ 0 & 0 & -278.6 & 50 & 0 \\ 0 & 0 & 0 & -288.6 & 0 \\ 0 & 0 & 100 & 0 & -268.6 \end{pmatrix} \begin{pmatrix} z_1 \\ z_2 \\ z_1 z_2 \\ z_1^2 \\ z_2^2 \end{pmatrix}$$

$$+ \begin{pmatrix} -1 & 0 & 0 & 0 & 0 \\ 0 & -1 & 0 & 0 & 0 \\ -1.12 & 7 & -2 & 0 & 0 \\ 14 & 0 & 0 & -2 & 0 \\ 0 & -2.24 & 0 & 0 & -2 \end{pmatrix} z\tilde{u} + \begin{pmatrix} 7 \\ -1.12 \\ 0 \\ 0 \\ 0 \end{pmatrix} \tilde{u}. \quad (24a)$$

$$y = (0 \ 1 \ 0 \ 0 \ 0)z. \quad (24b)$$

From the bilinear model in Eqs. (24a) and (24b) of the van de Vusse reactor it is relatively straight forward to derive the linear process model, second order as well as third order correction term. Using Eqs. (24a) and (24b) given above, the linear reactor model and second order Volterra model can be described by the following frequency domain Volterra kernel (Doyle *et al.* 1995), i.e. to achieve the control objective outlined above, controllers using Volterra models up to third order correction term were designed, along with the corresponding linear controller.

$$H_1(s) = c^T(sI_1 - A_{11})^{-1}b_1, \ H_2(s_1, s_2) = c^T((s_1 + s_2)I_2 - A)^{-1}N(s_1 - A)^{-1}b, \quad (25)$$

$$H_3(s_1, s_2, s_3) = c^T((s_1 + s_2 + s_3)I_3 - A)^{-1}N((s_1 + s_2)I_2 - A)^{-1}N(sI_1 - A)^{-1}b, \quad (26)$$

Combining Eqs. (23)-(26), we get the first order and second order Volterra model for van de Vusse reactor, i.e.

$$H_1(s) = \frac{188.384 - 1.12s}{s^2 + 278.6s + 19379.49}. \quad (27)$$

$$H_2(s_1, s_2) = \frac{(1.12s_1 - 188.384)(s_1 + s_2 + 144.3)(s_1 + s_2 + 288.6) - 350(s_1 + 134.4)(s_1 + s_2 + 428.6)}{(s_1 + 134.4)(s_1 + 144.3)(s_1 + s_2 + 134.4)(s_1 + s_2 + 144.3)(s_1 + s_2 + 288.6)}. \quad (28)$$

Third-order Volterra model equation

$$p_3(s_1, s_2, s_3) = \frac{350(s_1 + s_2 + 428.6)}{(s_1 + s_2 + 134.4)(s_1 + s_2 + 144.3)(s_1 + s_2 + 288.6)(s_1 + s_2 + s_3 + 134.4)(s_1 + 134.4)(s_1 + 144.3)}$$

$$+ \frac{350(s_1 + s_2 + 428.6)(s_1 + s_2 + 302.6)}{(s_1 + s_2 + 134.4)(s_1 + s_2 + 288.6)(s_1 + s_2 + s_3 + 144.3)(s_1 + s_2 + s_3 + 134.4)(s_1 + s_2 + s_3 + 288.6)(s_1 + 134.4)(s_1 + 144.3)}$$

$$+ \frac{98000}{(s_1+s_2+134.4)(s_1+s_2+s_3+134.4)(s_1+s_2+s_3+144.3)(s_1+s_2+s_3+288.6)(s_1+134.4)(s_1+144.3)}$$

$$+ \frac{188.384-1.12s}{(s_1+s_2+134.4)(s_1+s_2+s_3+134.4)(s_1+134.4)(s_1+144.3)}. \tag{29}$$

### 4.1 Open loop response of linearized van de Vusse reactor

The above mentioned Volterra models Eqs. (27)-(28) are realized as mentioned in *section 2.3* on Matlab/Simulink platform. Figure10 and Figure11 shows the open-loop response for concentration *B* stemming from three different system models to the step changes in the inlet flow rate of +15 and -20 respectively. The solid black line in Figure 10 and Figure 11 denotes the response of original nonlinear van de Vusse model. The red line in Figure 10 and Figure 11 denotes the response of first order Volterra model and the green and blue lines denote the second order and the third order Volterra model respectively. It is clear from Figure 10 &11 that third order model shows a closer response to the true nonlinear response in contrast to the other three responses.

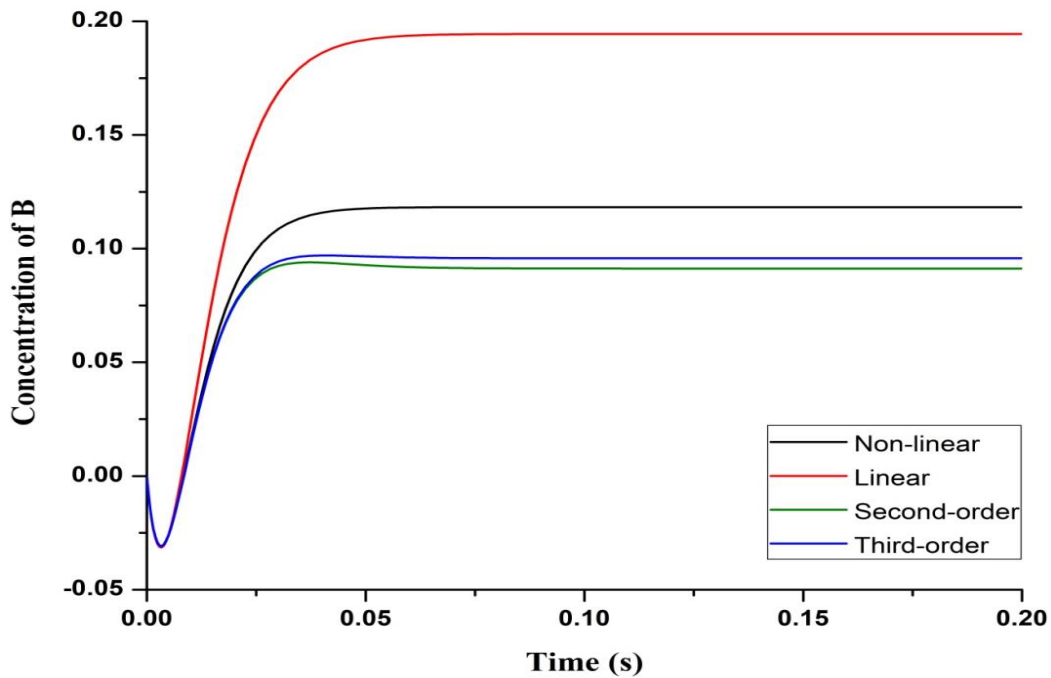

Fig.11: Response of Volterra based control model with positive step change

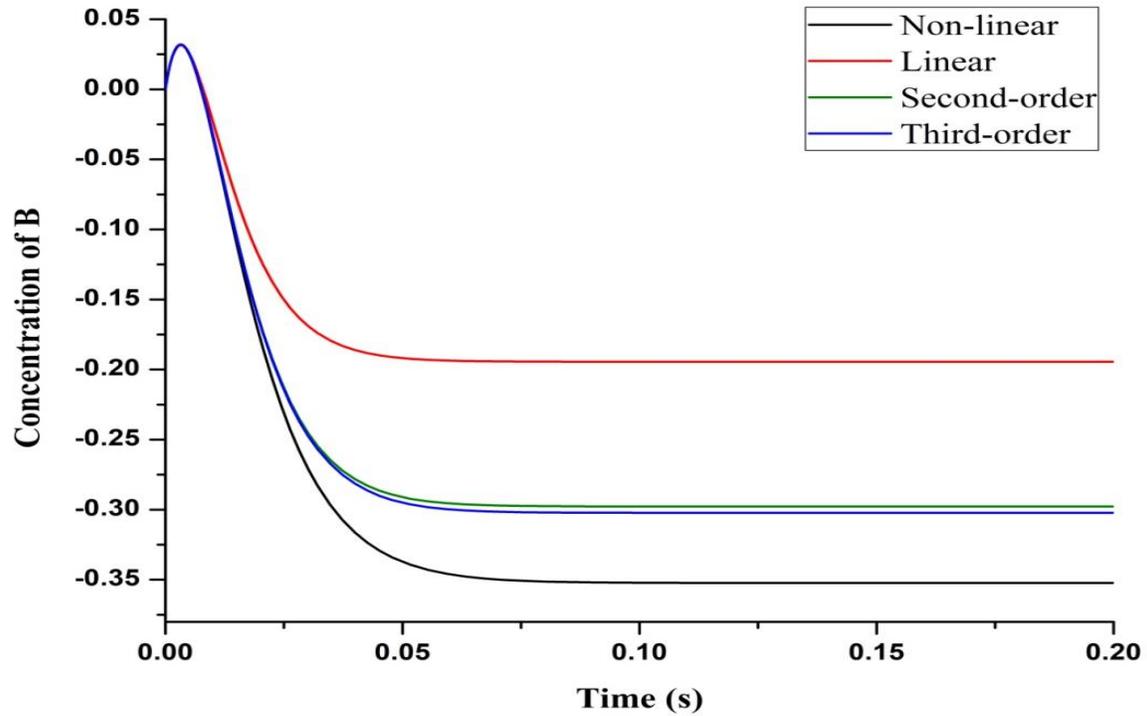

Fig.12: Response of Volterra based control model with negative step change

On a careful observation of the fig.11 and fig.12 it can be showed that as we increase the order of the derived Volterra models towards nth order the open loop trajectory goes near to the true nonlinear trajectory. A comparison of the trajectories depicted in Figure 10 and 11 is shown in Table (2). Table (2) shows closeness of the trajectories to the true nonlinear model. The values in the table (2) clearly indicates that Third order Volterra model is the closest to the nonlinear model with a difference of 18.85% and 14% respectively for positive and negative step change. Whereas, the linear model has a difference of 45% approximately which is quite away from the true trajectory and hence it is away from saying to be a good model.

**Table 1. Comparison of the open loop trajectories**

| | Steady-state values | | Difference in the steady-state from nonlinear trajectory (%) | |
|---|---|---|---|---|
| | For positive change in input flow rate (+20) | For negative change in input flow rate (-20) | +ve change | -ve change |
| Non-linear model | 0.1183 | -0.3523 | 0 | 0 |
| 1st order Volterra model (Linear) | 0.1944 | -0.1944 | 64.3 | 44.81 |
| 2nd order Volterra model | 0.09 | -0.2976 | 23.92 | 15.52 |
| 3rd order Volterra model | 0.096 | 0.3021 | 18.85 | 14 |

## 4.2 Controller synthesis and closed-loop results

Non-minimum phase systems pose a challenging task for control design since the exact system inverse is unstable. A van de vusse CSTR which is a non-minimum phase system can be represented as a partitioned nonlinear plant. To obtain an unconstrained nonlinear controller, all pass factorization method mentioned in *section 3* is applied to $H_1(s)$ which can be factored into an all pass and a minimum phase contribution (Morari and Zafiriou, 1989).

$$H_1(s) = H^A H^M, \qquad (30)$$

where $H^A$ is the all pass portion and $H^M$ is the minimum phase portion. Embedding Eq. (27) in to the above Eq. (30), we get

$$H_1^M = \frac{-1.12 - 188.348}{-(s^2 + 278.6s + 19379.49)} \text{ and } H_1^A = \frac{1.12 - 188.348}{-1.12 - 188.348}. \qquad (31)$$

Now utilizing Eq. (19) and embedding the above Eq. (31), we get the IMC-Volterra controller $V_{IMC}$ for the van de Vusse reactor.

$$V_{IMC}(s) = (H^M)^{-1}F.$$

$$= \frac{-(s^2+278.6s+19379.49)}{(-1.12-188.348)(\lambda s+1)}$$

$$= \frac{-(s^2+278.6s+19379.49)}{-0.0112s^2-1.12s-1.88348s-188.348}$$

$$= \frac{100s^2+27860s+1.93\times10^6}{1.12s^2-300.34s+18834.8}$$

$$= \frac{89.28s^2+24875s+1.729\times10^6}{s^2+268.167s-1.88348s+16816.78}. \tag{32}$$

The above equation (32) obtained is linear controller expression that is realizable. The rest of the Volterra models associated with the second and third order will act as the correction term to the IMC-Volterra controller.

Figure (13) shows the closed-loop result of the van de Vusse reactor control problem. Here, a positive step set point change is given in the reference value of concentration B. Figure (13) displays the ideal trajectory in comparison to the closed-loop trajectory generated by the third order Volterra model representation, second order Volterra representation and the linear model. In the linear model trajectory only the first order representation is considered with a linear controller. However, the second and third order trajectories consist of nonlinear model as well as their correction terms in the inner feedback loop of the controller.

Consider case 2, where a negative step set-point change in given into the reference value of the concentration of product B. Figure 14 depicts the closed-loop result of case 2. Figure 14 does not depict much of a difference into the trajectories of all the models in closed-loop settings. But, it is quite evident from the figure 14 that third-order Volterra model provides a much closer trajectory to the ideal trajectory.

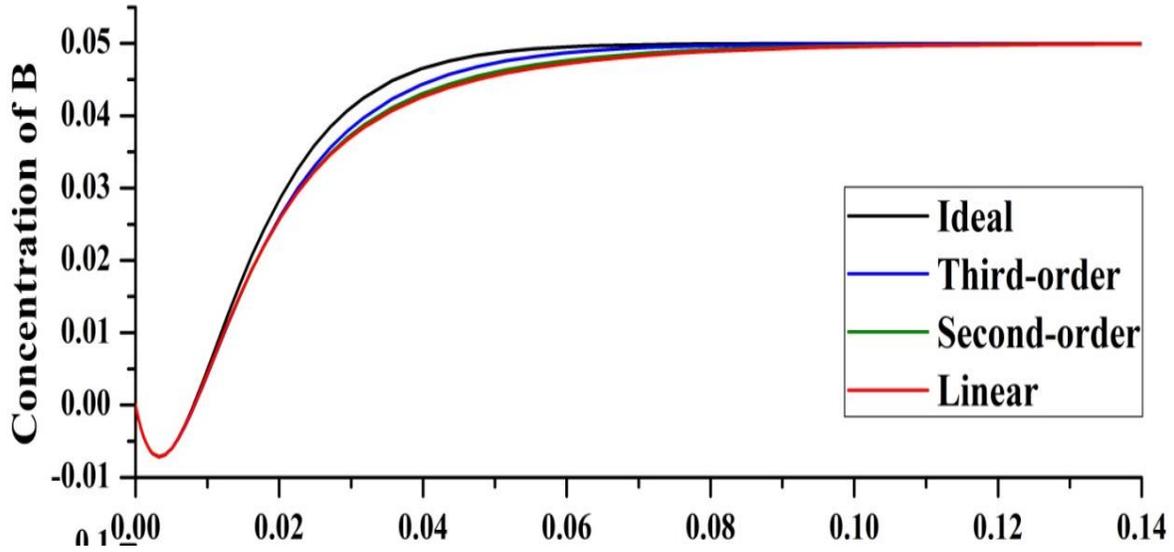

Fig.13: Closed-loop simulation for concentration of product B with positive step set-point change.

In order to evaluate the controller performance qualitatively, we adopt Integral Square Control Input (ISCI), i.e $ISCI = \int_0^t u^2(\tau)d\tau$. The less value of ISCI indicates less efforts required by the controller to achieve the desired output. To demonstrate the reduction in the control efforts through the closed-loop scheme associated with the third order Volterra model, ISCI values of all the closed-loop schemes are measured and noted in Table (3).

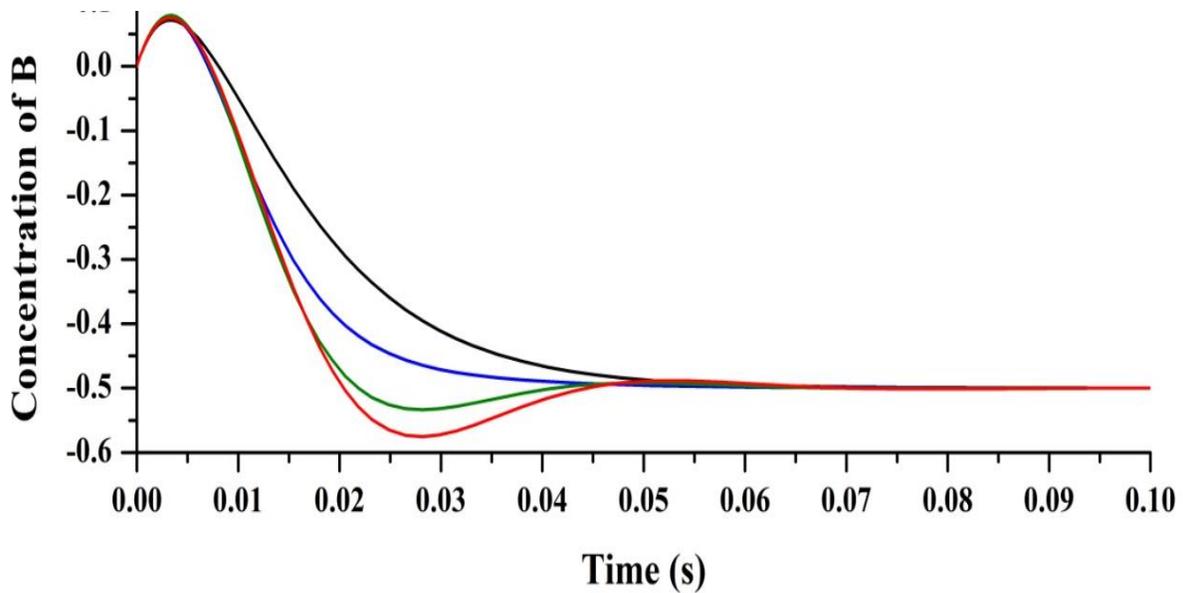

Fig.14: Closed-loop simulation for concentration of product B with negative step set-point change.

Table 2. ISCI index

| Method | Control Efforts in Inputs (ISCI) |
|---|---|
| Linear | 185 |
| Second order Volterra model | 183 |
| Third order Volterra model | 172.4 |

The values in Table (3) clearly shows the superiority of the third order IMC-Volterra controller associated with the third order Volterra model in contrast to the other two controllers. The ISCI value associated with the third order controller is significantly less than the linear controller by 7.02%. That is indicative of better performance from the proposed method.

## 5. Conclusion

In this paper, we achieve the Volterra model up to third order for van de Vusse reactor via introduction of Carleman linearization method. A general theory to arrive at the bilinear approximation for arbitrary nonlinear model is represented in this work. This paper will be useful for deriving the efficient $n^{th}$ order Volterra model for any arbitrary nonlinear systems. Moreover, results of control investigation by designing IMC controller based on Volterra model for van de Vusse reactor is presented. Despite the complexity of model, results of this paper show the superiority of the third order Volterra model which circumvent the complexity over the linear model and second order Volterra model and is attributed to its closeness to the true nonlinear model by capturing the dynamics of the system. Comparison of all the controllers is shown among which controller with third order model shows superior performance than the rest. Thus, IMC control using Volterra models in combination with method of realization for nonlinear systems is recommended. This method is also recommended for unstable non-minimum phase nonlinear systems.